\documentclass[11pt]{amsart} 
\usepackage{amssymb,hyperref}
\usepackage{graphicx}

\newcommand{\PP}{\mathbb P}

\newcommand{\CC}{{\mathbb C}}
\newcommand{\NN}{{\mathbb N}}
\newcommand{\FF}{{\mathbb F}}
\newcommand{\ZZ}{{\mathbb Z}}
\newcommand{\QQ}{{\mathbb Q}}
\newcommand{\RR}{{\mathbb R}}

\newcommand{\mA}{{\mathcal A}}

\newcommand{\mE}{{\mathcal E}}
\newcommand{\mo}{{\mathcal O}}

\newcommand{\Pic}{\mathrm{Pic}}
\newcommand{\NS}{\mathrm{NS}}
\newcommand{\MWL}{\mathrm{MWL}}
\newcommand{\Hod}{\mathrm{Hod}}

\newcommand{\End}{\mathrm{End}}

\theoremstyle{remark}

\theoremstyle{definition}

\theoremstyle{plain}

\numberwithin{equation}{section}
\begin{document}

\title{On families of K3 surfaces with real multiplication}

\author{Bert van Geemen}

\address{Dipartimento di Matematica, Universit\`a di Milano,
Via Saldini 50, I-20133 Milano, Italia}
\email{lambertus.vangeemen@unimi.it}

\author{Matthias Sch\"utt}
\address{Institut f\"ur Algebraische Geometrie, Leibniz Universit\"at
  Hannover, Welfengarten 1, 30167 Hannover, Germany}

    \address{Riemann Center for Geometry and Physics, Leibniz Universit\"at
  Hannover, Appelstrasse 2, 30167 Hannover, Germany}

\email{schuett@math.uni-hannover.de}

\begin{abstract}
We exhibit large families of K3 surfaces with real multiplication,
both abstractly using lattice theory, the Torelli theorem and the surjectivity of the period map,
as well as explicitly using dihedral covers and isogenies.
\end{abstract}

\date{November 1, 2023}

\maketitle

\section{Introduction}

Most of the theory of complex K3 surfaces is governed by the Hodge structure on the second cohomology group.
The symmetries of this Hodge structure lead
to the concepts of real multiplication (RM) and complex multiplication (CM).
While we already have a decent understanding of CM,
in particular for the existence problems which we consider here, by work of Taelman \cite{Taelman},
RM remains rather mysterious,
with only very few abstract constructions and even fewer concrete examples so far (cf.\ \ref{ss:RM})
and no analogue of Taelman's result.
The present paper aims to remedy this by developing new general methods
which can be used to construct families of K3 surfaces with RM (and with CM).
Here a family of K3 surfaces has RM (or CM) by a field $F$ if the very general member $X$ in the family has
$F=\End_{\Hod}(T_{X,\QQ})$. The K3 surfaces considered in this paper are all algebraic.

We obtain general results on the existence of families of maximal dimension of K3 surfaces with RM
by any totally real field of degree $m=2,5$ using
Taelman's results on K3 surfaces with CM, see Theorem \ref{thm:genRM}. A variation of his method
allows us to deal with the remaining cases $m=3,4,6,7$ in Theorems \ref{thm:cubic} and \ref{thm:RMmax367}.
In the case of real quadratic fields we give an alternative proof of the existence of such families,
with more explicit K3 surfaces and with an explicit description of the action of $F$ on the transcendental lattice,
in Theorem \ref{thm:quadRMmax}.

To get explicit examples of families of K3 surfaces with RM, our first
approach uses K3 surfaces $X$ with a purely non-symplectic automorphism $\sigma$
of order $m$ so that $X$ has CM with the cyclotomic field $\QQ(\zeta_m)$.
We deform the cyclic covering $X\rightarrow X/\sigma$ in such a way that the action of the totally real
subfield of $\QQ(\zeta_m)$ on the transcendental cohomology deforms with $X$.
These deformations are not Galois coverings, but their monodromy group is the dihedral group $D_m$ of order $2m$.
This provides a concrete implementation of the ideas that go into proof of Theorem \ref{thm:genRM}.
We obtain the following results.

\subsection{Theorem}
\label{thm2}
Let $\rho$ denote the rank of $\Pic(X)$ for a very general $X$ in a family,
so $d=\dim_\QQ T_{X,\QQ}=22-\rho$.
\begin{enumerate}
\item[(5)]
The  7-dimensional family of degree 2 K3 surfaces in \S \ref{ss:5big}
has $\rho= 2$ and RM by $\QQ(\sqrt 5) = \QQ(\zeta_5+\zeta_5^{-1})$.
\item[(7)]
The 3-dimensional family of elliptic K3 surfaces in \S \ref{ss:7}
has $\rho=4$ and RM by the cubic field $\QQ(\zeta_7+\zeta_7^{-1})$.
\item[(9)]
The 2-dimensional family of elliptic K3 surfaces in \S \ref{ss:9}
has $\rho=10$ and RM by the cubic field $\QQ(\zeta_9+\zeta_9^{-1})$.
\item[(11)]
The 2-dimensional family of elliptic K3 surfaces in \S \ref{ss:11}
has  $\rho= 2$ and RM by the degree five field $\QQ(\zeta_{11}+\zeta_{11}^{-1})$.
\end{enumerate}

Our second approach, in Sections \ref{s:isogeny}, \ref{s:higher},
exploits isogenies between elliptic K3 surfaces
to exhibit self-maps of K3 surfaces, a topic of independent interest.
Such a self-map is then shown to induce RM or CM.

\subsection{Theorem}     
\label{thm}
Let $\rho$ denote the rank of $\Pic(X)$ for a very general $X$ in a family,
so $d=\dim_\QQ T_{X,\QQ}=22-\rho$.
\begin{enumerate}
\item[(2)]
The 4-dimensional family of elliptic K3 surfaces in Proposition \ref{prop:2}
has $\rho= 10$ and RM by $\QQ(\sqrt 2)$.
\item[(3)]
The 3-dimensional family of elliptic K3 surfaces in Proposition \ref{prop:3}
has $\rho= 10$ and RM by $\QQ(\sqrt 3)$.
\end{enumerate}

In fact, the very same approach gives also new large families with CM
(see Propositions \ref{prop:2}, \ref{prop:3}) and subfamilies with larger CM fields (see \ref{ss:strata}).
There is also an isolated example with RM by the field $\QQ(\sqrt 7)$
(see  \ref{rem:5,7}).

\subsection{Comment}

After receiving a draft of the manuscript, Eva Bayer-Fluckiger  informed us that one could give 
  a complete description of the  real fields that give RM 
on maximal families of K3 surfaces,
extending  the ideas from Section \ref{s:RM-general}.

\subsection{Acknowledgements} We are indebted to Alice Garbagnati, Eva Bayer-Fluckiger and Simon Brandhorst
for very helpful discussions and comments.

\

\section{Hodge structures and moduli}
\label{s:basics}

\subsection{Hodge structures}

A complex K3 surface $X$ defines a simple polarized weight two Hodge structure
$$
T_X\,:=\, \Pic(X)^\perp\quad\subset\,H^2(X,\ZZ),\quad\mbox{let}\quad
T_{X,\QQ}\,:=\,T_X\otimes_\ZZ\QQ~.
$$
The polarized Hodge structure
$T_{X}$ has $\dim T^{2,0}_{X}=1$.
A $\QQ$-linear map $a:T_{X,\QQ}\rightarrow T_{X,\QQ}$
is an endomorphism of this Hodge structure if its complexification $a_\CC$ preserves the Hodge
decomposition $T_{X,\CC}=\oplus T_X^{p,q}$.
Zarhin \cite[Thm.\ 1.5.1]{Zarhin} showed that the $\QQ$-algebra of these endomorphisms
$$
F\,:=\,{\End}_{\Hod}(T_{X,\QQ})\,=\,
\big\{a\,\in\,\End(T_{X,\QQ}):\;a_\CC(T^{p,q}_{X})\,\subset\, T^{p,q}_{X}\,\big\}
$$
is either a totally real field or a CM field. 
In the latter case, we say that $X$ has complex multiplication (CM),
but the notion of real multiplication (RM) is usually reserved for those K3 surfaces where $F\neq \QQ$
and $F$ is totally real.


Since $F$ is a field, $T_{X,\QQ}$ is an $F$-vector space, and we write
$$
d\,=\,\dim_\QQ T_{X,\QQ},\qquad  m\,:=\,[F:\QQ],\qquad l\,:=\,d/m\,=\,\dim_FT_{X,\QQ}.
$$
Since $F\otimes_\QQ\CC=\oplus_\sigma \CC$, where $\sigma$ runs over the complex embeddings of $F$,
the action of $F$ on $T_X$ is diagonalizable.
The polarization on $T_{X,\QQ}$, which is induced by the intersection form on $H^2$, has the property
\begin{eqnarray}
\label{adjoint}
(ax, y) = (x, \bar a y) \qquad \forall \, x,y \in T_{X,\QQ}, \quad a \in F~,
\end{eqnarray}
where $\bar{a}$ denotes the complex conjugate of $a$.

\subsection{Remark}

Some authors (for example, \cite[Def.\ 9]{Taelman}),
define $X$ to have CM by a CM field $F$ only if moreover $l = \dim_FT_{X,\QQ}=1$; 
this setting is of special interest because it comes with the extra feature
that  $X$ can be defined
over some number field, like the singular K3 surfaces highlighted in \ref{ss:CM}. Our  focus is, however, on exhibiting 
\emph{families} of K3 surfaces with RM or CM, so we allow for $l>1$, in agreement with \cite{Zarhin}.

\subsection{Complex multiplication (CM)}
\label{ss:CM}
A notable example of a K3 surface $X$ with CM is a singular K3 surface,
i.e.\ $X$ has maximal Picard number $\rho(X)=20$. The CM arises here from the
positive definite quadratic form $Q$ on the
transcendental lattice $T_X$, in fact the CM field is $F=\QQ(\sqrt{-\det(Q)})$.
These K3 surfaces form isolated points in moduli.

Other examples of K3 surfaces with CM are easily exhibited
by considering surfaces with a purely-nonsymplectic automorphism $\sigma$ of order $n>2$, so
we require that $\sigma$ acts on a hlomorphic 2-form $\omega$ as multiplication by a primitive $n$-th root of unity
$\zeta_n$:
\[
\sigma^*\omega = \zeta_n\omega.
\]
The action of $\sigma^*$ on $T_{X,\QQ}$ gives this $\QQ$-vector space the structure of a
$\QQ(\zeta_n)$-vector space and thus $X$ has CM by a field $F$ containing $\QQ(\zeta_n)$.
We shall use this in \ref{s:RMellK3} to find examples of RM.

A CM field $F\subset \End_{\Hod}(T_X)$ defines an eigenspace decomposition of $T_{X,\CC}$ for the action
of $F$. There are $m$ eigenspaces, each of dimension $l$. Let
$V\subset T_{X,\CC}$ be the eigenspace containing $H^{2,0}(X)$.
The surjectivity of the period map and the Torelli theorem
imply that deforming the subspace $H^{2,0}$ in $V$ deforms $X$, with the given action of $F$,
in a family of dimension
\begin{eqnarray}
\label{eq:CM-dim}
l-1 = \,\dim_\CC \PP V \,=\,\dfrac{\mbox{rank}(T_X)}{[F:\QQ]}-1.
\end{eqnarray}

In the presence of a purely non-symplectic automorphism of order $n$, the degree of
$F$ over $\QQ$ is given by $\phi(n)$, where $\phi$ is the Euler totient function.
In particular, $X$ is isolated in moduli when rank$(T_X)=\phi(n)$.

 \subsection{Real multiplication (RM)}
 \label{ss:RM}
 
We shall now compare the CM setting with the RM case following   \cite{RM}.
If $X$ is a K3 surface with RM by a totally real field $F$,
then $T_X\otimes\RR$ admits a decomposition into eigenspaces for the $F$-action,
each of real dimension equal to $l = $ rank$(T_X)/[F:\QQ]$.
Consider the special eigenspace $T_\epsilon$
whose complexification contains the 2-form $\omega$, i.e.
$$
H^{2,0}(X)\,\subset\, T_\epsilon\otimes\CC.
$$
Then, by complex conjugation,  also $H^{0,2}(X)\subset T_\epsilon\otimes\CC$,
i.e.\ $T_\epsilon$ has signature $(2,l-2)$ whereas the other eigenspaces are negative definite.
The main overall restriction on RM structures is the following:

\subsection{Lemma \cite[Lemma 3.2]{RM}} \label{dimRM}

In the RM case, one has
$l = \dim_\RR T_\epsilon\geq 3$.

\medskip

Indeed, otherwise $T_\epsilon$ would be positive definite of dimension two
which implies that $T_{X,\QQ}$ admits extra endomorphisms,
just like for singular K3 surfaces, and in fact $X$ has CM
(compare the example in Proposition \ref{prop:RM5}).

\subsection{Moduli}
\label{ss:moduli}

As a consequence, a K3 surface with RM (by a field $F\neq \QQ$) has a transcendental lattice of rank at least $6 = 2\cdot 3$,
so the Picard number satisfies 
\begin{eqnarray}
\label{eq:rho}
\rho(X)\leq 16.
\end{eqnarray}
The deformation space of K3 surfaces with given $F$-action on $T_{X,\QQ}$ has dimension $l-2$
(see the proof of \cite[Lemma 3.2]{RM}):
\begin{eqnarray}
\label{eq:RM-dim}
l-2 \,=\, \dim_\RR T_\epsilon -2 \,=\, \frac{\mbox{rank}(T_X)}{[F:\QQ]}\,-\,2\,\geq\, 1~.
\end{eqnarray}
In particular, a K3 surface with RM is not isolated in moduli.

\subsection{Families of K3 surfaces with RM}
To conclude, if an algebraic K3 surface $X$ has RM by a field $F$
one has $d\leq 21$, hence $l\geq 3$ implies that $m=[F:\QQ]\leq 7$.
Besides the condition that $m\leq 7$,
we do not know of previous general results on the problem which totally real fields can be obtained
as ${\End}_{\Hod}(T_{X,\QQ})$ for a K3 surface $X$, besides the abstract ones with $[F:\QQ]=2$ and $l=3$
in \cite[Example 3.4]{RM} and the impressive explicit examples in
the papers \cite{ElsenhansJ1}, \cite{ElsenhansJ2} and \cite{ElsenhansJ3} for real quadratic fields and high Picard rank.

\subsection{Cyclotomic fields}

In this paper we consider in particular totally real subfields  $F\subset \QQ(\zeta_{n})$ of
the cyclotomic field of $n$-th roots of unity. This CM number field is a Galois
extension of $\QQ$ with Galois group the group of units $(\ZZ/n\ZZ)^\times$
in $\ZZ/n\ZZ$, it has degree $\phi(n)$ over $\QQ$,
where $\phi$ is again Euler's totient function.
Complex conjugation on $Gal(\QQ(\zeta_{n}))$ is given $-1\in (\ZZ/n\ZZ)^\times$
and thus the totally real subfields of $\QQ(\zeta_{n})$ correspond to the subgroups
$H<(\ZZ/n\ZZ)^\times$ with $-1\in H$;
in particular, any such field is contained
in the maximal totally real subfield $\QQ(\zeta_{n}+\zeta_{n}^{-1})$, of
degree $\phi(n)/2$ over $\QQ$.

\

\section{General results on RM}
\label{s:RM-general}

\subsection{Abstract deformations from CM to RM}
Given a K3 surface $X$ with CM by a field $E=\End_{\Hod}(T_{X,\QQ})$ of degree $m$ over $\QQ$,
the totally real subfield $F\subset E$, of degree $m/2$ over $\QQ$, also acts on $T_{X,\QQ}$.
In case $l=\dim_E (T_{X,\QQ})\geq 2$, one has $\dim_F(T_{X,\QQ})=2l\geq 4$, so the obstruction to RM in
Lemma \ref{dimRM} is not present and therefore $X$ is a member of a $2l-2$-dimensional
family of K3 surfaces with RM by $F$.

\subsection{Proposition}\label{prop:cyclo}
Let $X$ be a K3 surface with CM by the field $E$ and assume that
$l:=\dim_E (T_{X,\QQ})\geq 2$.

Then there exists a $2l-2$-dimensional family of K3 surfaces with real multiplication
by the totally real subfield $F$ of $E$.
These K3 surfaces are deformations
of $X$, and for the very general $X_\eta$ in this family there are
isometries $T_{X_\eta}\cong T_X$ and $\Pic(X_\eta)\cong \Pic(X)$.

\subsection{Proof} 

This follows from \ref{eq:RM-dim}, we provide some details (cf.\ \cite[Proof of Lemma 3.2]{RM}).
The inclusion $F\hookrightarrow \End_\Hod(T_{X,\QQ})$ induces a splitting of the
real vector space $T_{X,\RR}$ as a direct sum of
$[F:\QQ]$ eigenspaces for the action of $F_\RR:=F\otimes_\QQ\RR$, each of dimension $2l$.
Let
$$
V=T_\epsilon\subset T_{X,\RR}
$$
be the eigenspace with $H^{2,0}(X)\subset V_\CC$.
The intersection form on $H^2(X,\ZZ)$ induces a bilinear
form $(\cdot,\cdot)_V$ on $V$ of signature $(2+,(2l-2)-)$. The deformations $X_\eta$
of $X$ with $F\subset \End_\Hod(X_\eta)$ are parametrized by the one-dimensional subspaces
$\langle \omega_\eta\rangle$ of $V_\CC$ with $(\omega_\eta,\omega_\eta)_V=0$
and $(\omega_\eta,\overline{\omega_\eta})_V>0$. In this way one obtains a $2l-2$-dimensional
family of deformations of $X$.
Since $2l-2>0$ by assumption, the very general member in this family has $F=\End_\Hod(X_\eta)$, since overfields of $F$ give families with fewer moduli.
\qed

%

\subsection{Maximal families}
Recall that a K3 surface $X$ with RM by a field $F$ of degree $m(>1)$ has a transcendental lattice of rank
$d=l\cdot m$ for some $l\geq 3$ and $d\leq 21$ and then $X$ is a member of an $l-2$-dimensional family of K3 surfaces with RM by $F$. 
In case $m=2,5$ the maximal $l$ is $20/m=10,4$ respectively and Theorem \ref{thm:genRM} shows that,
for any totally real field of such a degree, a family of maximal dimension exists.

In case $m=3$ the maximal $l$ is $l=7$ whereas for $m=4$ it is $l=5$ and for  $m=6,7$ it is $l=3$. 
We will show in Theorems \ref{thm:cubic} and  \ref{thm:RMmax367}
that there exist K3 surfaces having RM with fields of these degrees with these values of $l$.
However, we were not able to find succinct conditions (or if there any non-trivial ones) on which totally real fields of these degrees and these $l$ occur as $\End_{\Hod}(X)$ for a K3 surface $X$.

\subsection{Theorem}\label{thm:genRM}
Let $F$ be a totally real field of degree $1<d\leq 5$.
Let $l=\lfloor\frac{20}{d}\rfloor$ for $d\neq 4$ resp.\ $l=4$ for $d=4$.
Then there is an $l-2$-dimensional family of K3 surfaces with RM by $F$.

\subsection{Proof} 
Given $F$, we embed it in a CM field $E$ as follows.
Let $K$ be any CM field of degree $l$ over $\QQ$ such that $K\cap F = \QQ$,
the intersection taken with respect to any embeddings $F, K \hookrightarrow\CC$.
(Note that $l>2$ is even, assuring the existence of such $K$; one can take $K$ to be composed of any imaginary quadratic field with almost any another field of complementary degree.)
Then the composite field  $E:=FK$ is a CM field of degree $dl\leq 20$ with $F\subset E$.
Taelman proved that there exists a  K3 surface $X$ with $E=\End_{\Hod}(T_{X,\QQ})$ and
$\dim_E(T_{X,\QQ})=1$ (see \cite[Thm.\ 4]{Taelman}).
Then $\dim_F(T_{X,\QQ})=l\geq 4$, and as $l\geq 3$,
this guarantees the existence of the family of the K3 surfaces that are deformations of $X$, with RM by $F$.
\qed

\

The examples we construct below use the following proposition.  It is a converse for the results discussed in the previous section.

\subsection{Proposition} Let $V$ be a $\QQ$ vector space with a non-degenerate bilinear form $(\cdot,\cdot)$ such that the quadratic form it defines on $V_\RR$
has signature $(2,n-2)$.
Let $F$, $F\neq\QQ$, be a totally real number field and assume  that
$V$ also has the structure of an $F$-vector space such that the following two conditions hold true:
\begin{itemize}
\item
there is an eigenspace $V_\epsilon\subset V_\RR$ for all $a\in F$ on which the signature
of the restriction of the bilinear form is $(2,e)$ for some $e$;
\item the adjoint property 
 $(ax,y)=(x,ay)$ holds for all $a\in F$ and all $x,y \in V$.
 \end{itemize}

Then a positive definite oriented 2-plane in $V_\epsilon$ defines a simple Hodge structure of K3 type on $V$  with $F\subset \End_\Hod(V)$.
If $\dim_FV\geq 3$ then $F=\End_\Hod(V)$ for a very general $V$, so $V$ has RM by $F$.

\subsection{Proof} Let $\omega\in V_{\epsilon,\CC}$ be an eigenvector for the rotation by $\pi/2$ in the 2-plane such that the orientation is given by
$\omega+\bar\omega, (1/i)(\omega-\bar\omega)$. The Hodge structure on $V$ is defined by:
$$
V^{2,0}\,=\,\CC\omega,\quad V^{0,2}\,=\,\CC\bar\omega,\quad V^{1,1}\,=\,\langle\omega,\bar\omega\rangle^\perp\subset V_\CC.
$$
Then $(V,(\cdot,\cdot))$ is a polarized weight two Hodge structure of K3 type. To see that  $F\subset \End_\Hod(V)$ notice that $a V^{2,0}=V^{2,0}$, $aV^{0,2}=V^{0,2}$ since their bases lie in an eigenspace of $a$. Next we use the adjoint property: for $x\in V^{1,1}$ one has
$$
(ax,\omega)\,=\,(x,a\omega)\,=\,(x,\epsilon(a)\omega)\,=\,\epsilon(a)(x,\omega)\,=\,0,
$$
and similarly for $\bar\omega$.
Hence $ax\in V^{1,1}$ and we have  $F\subset \End_\Hod(V)$. That for a very general 2-plane one has equality follows by considering the 
Mumford Tate group of the Hodge structure $V$ as in  \cite{RM}. 
\qed

\subsection{Maximal families of K3 surfaces with RM by quadratic fields}
For any real quadratic field $F$, we can directly show the existence of algebraic K3 surfaces with a genus one fibration that have RM by $F$.
Such a K3 surface $X$ must have $d=\dim T_{X,\QQ}\leq 20$. Thus the maximal dimension
of a family of such K3 surfaces is $l-2=(20/2)-2=8$.
The existence of families of this dimension is a consequence of Theorem \ref{thm:genRM}.
Here we provide an alternative proof which provides families
with non-trivial geometrical information on the K3 surfaces $X$
and an explicit description of the action of $F$ on $T_X$.
We do not know of other families of maximal dimension, but their existence is quite likely.

\subsection{Theorem} \label{thm:quadRMmax}
\label{lem:d}
For any squarefree $d>0$ and any $r>0$, there is an 8-dimensional family
of K3 surfaces with a genus one fibration
such that the very general member $X$ has
$$
\Pic(X)\,\cong\,U(r)\,:=\,\left(\ZZ^2,\,\left(\begin{smallmatrix} 0&r\\r&0\end{smallmatrix}\right)\right)
$$
and has RM by $\QQ(\sqrt d)$.

\subsection{Proof} Consider a K3 surface with $\Pic(X)=U(r)$,
then $T_X = U\oplus U(r) \oplus E_8^2$, where $E_8$ is the unique unimodular negative definite lattice of rank $8$.
The linear system of an isotropic vector in $\Pic(X)$ endows $X$
with a genus one fibration which very generally only admits multisections of 
degree divisible by $r$.
The theorem thus follows as an application of the surjectivity of the period map and the Torelli theorem
once we endow $T_{X,\QQ}$ with a suitable action by $\QQ(\sqrt{d})$.
This can be achieved as follows. 
On any $U(r) \; (r\in\ZZ, r\neq 0)$, we define an endomorphism
$$
\tau: \; U(r) \,\longrightarrow\, U(r),\qquad
(u,v) \,\longmapsto\,  (dv, u),
$$
such that $\tau^2 = d$.
Notice that, for the diagonal action, the eigenspace
$$(
U_\RR+ U(r)_\RR)^{\tau=\sqrt d}\; =\; \RR (1, \sqrt d, 0, 0)\,  \oplus  \,\RR (0,0,1,\sqrt d)
$$
has indeed signature $(2,0)$, as required by \ref{ss:RM}.
On $E_8^2$ we consider the same map, $(u,v)\mapsto (dv,u)$ now with $u,v\in E_8$.
The diagonal action of these maps $M:T_{X,\QQ}\rightarrow T_{X,\QQ}$
satisfies $M^2=d$ and thus defines an action of $\QQ(\sqrt{d})$ which satisfies $q(Mx,y)=q(x,My)$
where $q(\cdot,\cdot)$ is the polarization on $T_{X,\QQ}$. The family of K3 surfaces is the one whose periods
lie in an eigenspace of $M$ in $T_{X,\CC}$.
\qed

\subsection{Remark}
Alternatively one can use the following method for $E_8^2$.
To find real
multiplication structures by all $\sqrt{d}$ on a given lattice $L$,
it suffices to assume that Aut$(L)$  contains 4 anti-commuting involutions
$g_1,\hdots,g_4$.
Then (in $\End(L)$)
\[
(a_1g_1+\hdots+a_4g_4)^2 \,=\, a_1^2+\hdots+a_4^2 \qquad \forall \, a_i\in \ZZ,
\]
so the endomorphism  $M=a_1g_1+\hdots+a_4g_4$, which satisfies \ref{adjoint}, endows $L_\QQ$ with
the structure of $\QQ(\sqrt{\sum a_i^2})$ vector space.
Using Lagrange's four-square theorem, this gives RM on $L$ by any real $\QQ(\sqrt d)$.
To conclude, one verifies that the required anti-commuting involutions can be found in the Weyl group $W(E_8)$.
We omit the details.

\

\subsection{Maximal families of K3 surfaces with RM by fields of degree $3, 4, 6, 7$}

For a K3 surface with RM by the field $F$, after identifying $T_\QQ=F^l$,
the property \ref{adjoint} is equivalent to the property that the intersection form on $H^2(X,\QQ)$,
restricted to $T_{X,\QQ}$, is given by
$$
(x,y)\,=\,Trace_{F/\QQ}({}^tx\Delta y),\qquad x,y\,\in\,F^l~,
$$
for some $l\times l$ matrix $\Delta$ with coefficients in $F$. It is easy to see that any such  bilinear form  has the property  \ref{adjoint} and for the converse one can argue as in the proof of
\cite[Proposition 9.2.3]{Birkenhake_L}.

To find the examples, we will restrict ourselves to the case that 
$T_{X,\QQ}$ is $\QQ$-isometric to $(\QQ^d, \mathrm I_{p,q})$, where $\mathrm I_{p,q}$ is the diagonal matrix with
$p$ diagonal coefficients equal to $+1$ and others equal to $-1$.
Besides  `well-known' fields given in Theorem \ref{thm:RMmax367},
we found many more examples using the same methods. However, only for cyclic cubic fields did we find a general lattice criterion
for the construction of K3 surfaces with RM. 
We discuss these fields first to explain our approach.

\subsection{Theorem}
\label{thm:cubic}
Let $F$ be a  totally real cyclic cubic field with class number one.
Then there is a 7-dimensional family of K3 surfaces with RM by $F$.

\subsection{Proof}
As we have seen, it suffices to find an action of $F$ on $T_{X,\QQ}$,
satisfying \eqref{adjoint} and the signature condition from \ref{ss:RM}, for some algebraic K3 surface $X$ with Picard rank one.
As $\Pic(X)=\ZZ h$ with $h^2=e$ for an even positive integer $e$,
one finds that $T_X\cong \ZZ v\oplus U^2\oplus E_8(-1)^2$ with $v^2=-e$.
In case $e=k^2$ is a square, $(1/k)v\in T_{X,\QQ}$ has square $-1$.
Hence the lattice generated by $(1/k)v$ and $U^2\oplus E_8(-1)^2$ is odd and unimodular, and thus,
since $\dim T_{X,\QQ}=21$, we have isometries
\begin{eqnarray}
\label{eq:TT}
T_{X,\QQ}&\cong&\langle 1\rangle^2\,\oplus\, \langle-1\rangle^{19} \nonumber
\\
&\cong&(\langle 1\rangle\,\oplus\, \langle-1\rangle^{2})^2\,\oplus\,(\langle-1\rangle^3)^5~.
\end{eqnarray}
Let $\mo_F$  be the ring of integers of $F$, it is a free $\ZZ$-module of rank $3=[F:\QQ]$.
For $y\in F$, the values of $Trace_{F/\QQ}(xy)$ are in $\ZZ$ for all $x\in \mo_F$
iff $y\in {\mathfrak d}^{-1}$ where ${\mathfrak d}\subset\mo_F$ is an ideal
called the different. So the dual lattice of $\mo_F$
with respect to the trace form
 is the fractional ideal $\mathfrak d^{-1}\subset F$.
By assumption, ${\mathfrak d}$ is a principal ideal and we choose a generator $\delta\in\mathcal O_F$.
Then the bilinear form $b_\delta$ on $\mo_F$ defined by
$$
b_{\delta}:\,\mo_F \times \mo_F\,\longrightarrow\,\ZZ,\qquad
b_{\delta}(x,y)\,:=\,Trace_{F/\QQ}(\delta^{-1} xy)
$$
is unimodular and thus there is an isometry $(F,b_\delta)\cong (\QQ^3, \mathrm I_{p,q})$
for certain $p+q=3$. For any unit $u\in \mo_F$, also $u\delta$ is a generator and we use this to find bilinear forms
$b_{u\delta}$ with the correct signatures to match \eqref{eq:TT}.

Recall that $F\otimes_\QQ\RR\cong\RR^3$, where $x\otimes\lambda$ maps to $\lambda (\sigma_1(x),\ldots,\sigma_3(x))$
where the $\sigma_i$ are the three embeddings $F\hookrightarrow\RR$.
As 
$$
Tr_{F/\QQ}(\delta^{-1}x^2)\,=\,\sum\sigma_i(\delta^{-1})\sigma_i(x)^2,
$$  
the signature of the $\RR$-linear extension of this quadratic form is determined by the signs of the $\sigma_i(\delta)$.
It follows from \cite{ArmitageF} (recently extended to higher degree fields in \cite[Cor.\ 4.3.5]{BVV})
 that $F$ has full unit signature rank.
That is, for any $p,q$ with $p+q=3$, there is a $u\in\mathcal O_F^\times$ such that
$\sigma_i(u\delta)$ assumes $p$ positive and $q$ negative values for $i=1,\ldots,3$.
In particular, we may assume that $b_{\delta}$ is negative definite and that $b_{u\delta}$ has
signature $ (1+,2-)$ for a certain unit $u$.

The isometry $T_{X,\QQ}\cong (F,b_{u\delta})^2\oplus(F,b_{\delta})^5$ and the diagonal action of $F$
 on the right hand side gives the desired action of $F$ on $T_{X,\QQ}$.
\qed

\subsection{Remark}
The crucial property that $\mathfrak d$ is a principal ideal also holds
generally whenever there is a single element $\alpha\in F$ such that $\mathcal O_F = \ZZ[\alpha]$.
For these cases the above proof applies
whenever the totally real cyclic cubic field $F$ has odd class number since it allows one to use \cite[Cor.\ 4.3.5]{BVV}
again.

\subsection{Example}
In the totally real cubic subfield $F=\QQ(\alpha)$ of  $\QQ(\zeta_7)$ 
with $\alpha=\zeta+\zeta^{-1}$, one can take 
$\delta=-\alpha^2 - 3\alpha -4$ totally negative, so $(F,b_{\delta})\cong (\QQ^3,\mathrm I_{0,3})$
and 
$u\delta = 2\alpha^2 - \alpha - 6$ such that 
$(F,b_{u\delta})\cong (\QQ^3,\mathrm I_{1,2})$.

\medskip

Using the same approach as above, we now cover real fields of the remaining degrees:

\subsection{Theorem} \label{thm:RMmax367}
For $m=4,6,7$ and $l=5,3,3$ respectively, there exist $(l-2)$-dimensional families
of K3 surfaces such that the very general member $X$ has RM by the field $F_m$
where $F_m$ is defined by the polynomial $f_m$ or $g_m$ in Table \ref{tab}.

\subsection{Proof}
We start by setting up the K3 surfaces in question
by specifying their Picard and transcendental lattices.

For $m=7$, we take exactly the same lattices as in the proof of Theorem \ref{thm:cubic}.

For $m=4$, let $r\in\NN$ and $\Pic(X)=U(r)$ as in Theorem \ref{thm:quadRMmax}.
Then, as a quadratic space,
\[
T_{X,\QQ} \cong U \oplus \langle 4r^2\rangle \oplus \langle -4r^2\rangle \oplus E_8^2 \cong \langle 1\rangle^2\oplus\langle-1\rangle^{18}~.
\]

For $m=6$, let $r\in\NN$ and $\Pic(X) = U \oplus \langle -4r^2\rangle^2$.
Then $T_X \cong \langle 4r^2\rangle^2 \oplus E_8^2$, so as a quadratic space
\[
T_{X,\QQ} \cong \langle 1\rangle^2 \oplus E_8^2 \cong  \langle 1\rangle^2\oplus\langle-1\rangle^{16}~.
\]
In summary, each case has $\dim T_{X,\QQ}=lm$ and thus
\begin{eqnarray*}
\label{eq:TT'}
T_{X,\QQ} &\cong& (\langle 1\rangle\,\oplus\, \langle-1\rangle^{m-1})^2\,\oplus\,(\langle-1\rangle^m)^{l-2}~.
\end{eqnarray*}
As in the proof of Theorem \ref{thm:cubic}, 
it remains to find a totally negative generator $\delta$ of $\mathfrak d$ and a unit $u\in\mathcal O_F^\times$
such that $(F,b_{u\delta})$ is hyperbolic.

This can be achieved with the help of  Magma \cite{Magma} as follows.
The table below gives a defining polynomial $f_m$ or $g_m$ of the field $F_m$.
Let $\alpha\in F_m$ be a root of $f_m$ resp.\ $g_m$;
then we write our choices of $\delta$ and $u\delta$ as linear combinations of the $\alpha^i$.

\begin{table}[ht!]
$$
\begin{array}{l}
 \begin{array}{rcl}
f_{4}&=& x^4 + x^3 - 6x^2 - x + 1,\\
\delta&=& -14\alpha^3 - 19\alpha^2 + 76\alpha + 46,\\
u\delta&=&(3\alpha^3 - 2\alpha^2 - 26\alpha - 5)/2,
\end{array}\\

\\
 \begin{array}{rcl}
g_{4}&=& x^4 - 6x^2 + 4,\\
\delta&=& 3\alpha^3 - 2\alpha^2 - 14\alpha + 16,\\
u\delta&=&\alpha^3 - 2\alpha^2 - 8\alpha + 16,
\end{array}\\

\\
 \begin{array}{rcl}
f_{6}&=& x^6 + x^5 - 5x^4 - 4x^3 + 6x^2 + 3x - 1,\\
\delta&=& 10\alpha^5 + 4\alpha^4 - 55\alpha^3 - 7\alpha^2 + 72\alpha - 21,\\
u\delta&=&
-2\alpha^5 + 7\alpha^4 + 11\alpha^3 - 22\alpha^2 - 4\alpha - 1,
\end{array}\\
\\
 \begin{array}{rcl}
g_7&=& x^7 + x^6 - 18x^5 - 35x^4 + 38x^3 + 104x^2 + 7x - 49,\\
\delta&=&  
(48274\alpha^6 - 25217\alpha^5 - 830561\alpha^4 - 425142\alpha^3 + 2481943\alpha^2 + 1241898\alpha -
    1553433)/7\\
    u\delta&=&  
(138\alpha^6 - 582\alpha^5 - 1418\alpha^4 + 6270\alpha^3 + 1455\alpha^2 - 15145\alpha + 7749)/7
\end{array}
\\

\\
 \begin{array}{rcl}
f_7&=& x^7 - 2x^6 - 5x^5 + 9x^4 + 7x^3 - 10x^2 - 2x + 1,\\
\delta&=&  
2\alpha^6 - 3\alpha^5 - 13\alpha^4 + 7\alpha^3 + 15\alpha^2 - 3\alpha - 6,  \\
u\delta&=&  
4\alpha^6 - 13\alpha^5 - 2\alpha^4 + 33\alpha^3 - 11\alpha^2 - 20\alpha - 6
\end{array}
\end{array}
$$
\caption{Totally real fields with suitable generators of the different ideal}
\label{tab}
\end{table}

The fields defined by $f_3,f_6$ are the totally real subfields of $\QQ(\zeta_7)$ and $\QQ(\zeta_{13})$ respectively,
with $\alpha=\zeta+\zeta^{-1}$. 

The field defined by $f_4$ is the degree 4 totally real subfield of $\QQ(\zeta_{17})$,
with $\alpha=\zeta+\zeta^4+\zeta^{-4}+\zeta^{-1}$, thus cyclic over $\QQ$,
while $g_4$ defines the biquadratic field $\QQ(\sqrt 2, \sqrt 5) = \QQ(\alpha)$ for  $\alpha  = \zeta^3+\zeta^{13}+\zeta^{-13}+\zeta^{-3}$ 
inside $\QQ(\zeta_{40})$.

The field defined by $f_7$ is not a Galois extension,
but the field defined by $g_7$ is Galois as it is the degree 7 totally real subfield of $\QQ(\zeta_{43})$,
with $\alpha=\zeta+\zeta^6+\zeta^7+\zeta^{-6}+\zeta^{-7}+\zeta^{-1}$.

In each case, we have by the choice of $\delta, u$ that
\[
T_{X_\QQ} \cong (F,b_{u\delta})^2 \oplus (F,b_{\delta})^{l-2}.
\]
Endowing this with the natural diagonal $F$-action gives the theorem.
\qed

\subsection{Remark}
In Sections \ref{s:ex}, \ref{s:isogeny}, \ref{s:higher}
we will exhibit explicit families of RM K3 surfaces using jacobian elliptic fibrations.
However, the dimensions of the families will sometimes be smaller then the maximum allowed for
by the degree of the field and the Picard number of the general member.



\section{Dickson polynomials and $D_n$-type covers $\PP^1\rightarrow\PP^1$}
\label{s:Dickson}


\subsection{The dihedral group $D_n$}
Let $D_n$ be the dihedral group of order $2n$.
We denote by $\tau, \sigma\in D_n$ an element of order two, $n$ respectively so that
$$
D_n\,=\,\langle\sigma,\tau\,\; |\;\, \sigma^n\,=\,\tau^2\,=\,1,\quad
\sigma\tau\,=\,\tau\sigma^{-1}\rangle~.
$$
We refer to \cite{CataneseP} for  general results on $D_n$-covers of algebraic varieties.
The $D_n$-action on certain varieties allows
us to obtain K3 surfaces with RM by the totally real subfield of $\QQ(\zeta_n)$.
For this we need the deformations of cyclic covers provided by the Dickson polynomials.

\subsection{Dickson polynomials and $D_n$-type covers of $\PP^1$}
The Dickson polynomial of degree $n$ with parameter $a$ is the (unique) degree $n$ polynomial
$p_{n,a}(x)\in \ZZ[a][x]$
satisfying, in the Laurent ring $\ZZ[a][v,v^{-1}]$,
$$
p_{n,a}(v+a/v)\,=\,v^n\,+\,(a/v)^n~.
$$
In particular, for $a=0$ we get $p_{n,0}=x^n$, so $p_{n,0}:\PP^1_x\rightarrow\PP^1_u$, $u=x^n$, is a
cyclic cover.
One easily verifies that
$$
p_{n,a^2}(ax)\,=\,a^np_{n,1}(x)~.
$$

The following Dickson polynomials will be used in this paper:
$$
\begin{array}{rcl}
p_{3,a}&=&x^3-3ax,\\
p_{5,a}&=&x^5-5ax^3+5a^2x,\\
p_{7,a}&=&x^7 - 7ax^5 + 14a^2x^3 - 7a^3x,\\
p_{9,a}&=&x^9 - 9ax^7 + 27a^2x^5 - 30a^3x^3 + 9a^4x,\\
p_{11,a}&=&x^{11} - 11ax^9 + 44a^2x^7 - 77a^3x^5 + 55a^4x^3 - 11a^5x~.\\
\end{array}
$$

\subsection{Lemma} \label{DicksonDn}
For any $n>2$ and for any non-zero $a\in\CC$,
the map defined by a degree $n$ Dickson polynomial $p_{n,a}$,
$$
f:\,\PP^1_x\longrightarrow\,\PP^1_u,\qquad x\,\longmapsto\, u\,:=\,p_{n,a}(x)~,
$$
is a degree $n$ covering with monodromy group $D_n$.  This covering is totally ramified over $\infty\in\PP^1_u$.

\subsection{Proof}

For a non-zero $a$ we define an action of $D_n$ on $\PP^1_v$ by
$$
\sigma,\tau:\,\PP^1_v\,\longrightarrow\,\PP^1_v,\qquad
\sigma:\,v\,\longmapsto\,\zeta_nv,\qquad \tau:\,v\,\longmapsto\,a/v~.
$$
Consider the composition
$$
\tilde{f}:\,\PP^1_v\,\longrightarrow\,\PP^1_x\,\longrightarrow\, \PP^1_u,\qquad u\,:=\,v^n\,+\,(a/v)^n,
$$
where the first map has degree two and is given by $x=v+a/v$. Then one finds that $\PP^1_x=\PP^1_v/\tau$ and  $\PP^1_u=\PP^1_v/D_n$.
Therefore the degree $n$ map $f:\PP^1_x\rightarrow\PP^1_u$
has monodromy group $D_n$.
\qed

\

\subsection{Real multiplication on elliptic K3 surfaces}\label{s:RMellK3}

We apply these deformations to cyclic covers of prime degree of jacobian elliptic surfaces.

Let $\overline\mE\rightarrow \PP^1_t$ be an elliptic surface  and assume that its Weierstrass model is
defined by a minimal Weierstrass equation (with $\alpha,\beta\in\CC[u]$):
\begin{eqnarray}
\label{eq:Weier_mEx}
\overline{\mE}:\quad Y^2\,=\,X^3\,+\,\alpha(u)X\,+\,\beta(u)~.
\end{eqnarray}
Then $\overline\mE$ is rational if $\deg(\alpha)\leq 4$ and $\deg(\beta)\leq 6$.
(cf.\ \cite[Prop.\ 5.51]{SchuettS}). Similarly, $\overline\mE$ is a K3 surface if it is not rational,
$\deg(\alpha)\leq 8$ and $\deg(\beta)\leq 12$ and the fibration is relatively minimal.
In particular, if for $n>1$ the equation
\begin{eqnarray}
\label{eq:Weier_mE}
\mE:\quad Y^2\,=\,X^3\,+\,\alpha(x^n)X\,+\,\beta(x^n)~.
\end{eqnarray}
defines a K3 surface, then (\ref{eq:Weier_mEx} defines a rational surface which is the quotient
$\mE/\sigma_0$ (with quotient map that sends $x\mapsto u=x^n$)
by the purely non-symplectic automorphism $\sigma_0$ of order $n$ given by
\begin{eqnarray}
\label{eq:sigma0}
\sigma_0:\,\mE\,\longrightarrow\,\mE,\qquad \sigma_0(X,Y,x)\,=\,(X,Y,\zeta_nx)~.
\end{eqnarray}
That $\sigma_0$ is non-symplectic can also be seen by computing $\sigma_0^*$
of the regular 2-form $dX\wedge dx/Y$ on $\mE$.

\subsection{Proposition}
\label{prop:RM-ell}

Let $n\in\{5,7,11\}$ and $a\in\CC^\times$. Assume that $\mE$, defined by (\ref{eq:Weier_mE}) is a K3 surface.
Then the deformation $\mE_a$ of  $\mE$
defined by the Weierstrass equation
\begin{eqnarray}
\label{eq:E_a}
\mE_a:\quad Y^2\,=\,X^3\,+\,\alpha(p_{n,a}(x))X\,+\,\beta(p_{n,a}(x)),
\end{eqnarray}
is again an elliptic K3 surface and the totally real subfield $F=\QQ(\zeta_n+\zeta^{-1}_n)$
of $\QQ(\zeta_n)$ acts by Hodge endomorphisms on $T_{\mE_a}$, so $F\subset \End_{\Hod}(T_{\mE_a})$.

\subsection{Proof}
To see this, let $\tilde{\mE}_a$ be the (relatively minimal) elliptic surface obtained as
the pull-back of $\mE_a$ along the double cover $\tilde{f}:\PP^1_v\rightarrow\PP^1_x$ defined by
$x:=v+a/v$. It has a Weierstrass equation 
$$
\tilde\mE_a:\quad Y^2=X^3+\alpha(v^n+(a/v)^n)X+\beta(v^n+(a/v)^n)~.
$$
Then $D_n$ acts on the Weierstrass model via its action on $\PP^1_v$, i.e.\ by
$$\sigma(X,Y,v)=(X,Y,\zeta_nv), \;\;\; \tau(X,Y,v)=(X,Y,a/v).
$$
This action extends to $\tilde\mE_a$, and $\tilde\mE_a/D_n$ is birational to the rational surface  $\overline\mE$.

Consider the Hodge substructure $T_{\mE_a,\QQ}\subset H^2(\mE_a,\QQ)$ defined by the transcendental lattice of
the K3 surface $\mE_a$, it is simple and has $\dim_\CC T_{\mE_a}^{2,0}=1$.
After pull-back to a desingularization
of the base change and push-foward along blow-downs to obtain the relative minimal model $\tilde{\mE}_a$,
one obtains
a Hodge substructure $T_a\subset H^2(\tilde{\mE}_a,\QQ)$ with an isomorphism of Hodge structures
$T_a\stackrel{\cong}{\longrightarrow}T_{\mE_a,\QQ}$.
It suffices to show that $F\subset \End_{\Hod}(T_a)$.

Since the rational map $\tilde{\mE}_a\rightarrow {\mE}_a$ is birational to the quotient by $\tau$,
one has $T_a\subset H^2(\tilde{\mE}_a,\QQ)^{\tau^*}$, the subspace of $\tau^*$-invariants, and
$H^{2,0}(\tilde{\mE}_a)^{\tau^*}\cong H^{2,0}(\mE_a)$.
Since ${\mE}_a$ is a K3 surface, we see that $H^{2,0}(\tilde{\mE}_a)^{\tau^*}$ is one dimensional
and hence $T_a$ is the unique simple Hodge substructure of $H^2(\tilde{\mE}_a,\QQ)^{\tau^*}$
with non-zero $(2,0)$-component.

Since $\sigma\tau\,=\,\tau\sigma^{-1}$, the endomorphisms $\sigma^*+(\sigma^{-1})^*$
and $\tau^*$ of $H^2(\tilde{\mE}_a,\QQ)$ commute. Thus $\sigma^*+(\sigma^{-1})^*$ defines an endomorphism of
the subspace of $\tau^*$-invariants $H^2(\tilde{\mE}_a,\QQ)^{\tau^*}$. This endomorphism maps $T_a$ into itself,
since the automorphisms $\sigma,\sigma^{-1}$ preserve the Hodge structure, and by the unicity of $T_a$. In particular, $\sigma^*+(\sigma^{-1})^*\in\End_\Hod(T_a)$.

Since $(\sigma^*)^n$ is the identity on $H^2(\tilde{\mE}_a,\QQ)$ and $n$ is prime, the subalgebra of $\End(H^2(\tilde{\mE}_a,\QQ))$
it generates is a quotient of $\QQ[T]/(T^n-1)\cong\QQ\times\QQ(\zeta_n)$.
The subalgebra generated by $\sigma^*+(\sigma^{-1})^*$ is therefore a quotient of $\QQ\times F$ and $F\not\cong\QQ$ since $n>3$.
To show that $T_a$ is an $F$-vector space it suffices to show that
$\sigma^*+(\sigma^{-1})^*$ acting on $T_a$ has no eigenvalue $\lambda\in \QQ$.
An eigenspace of $\lambda\in\QQ$ is a Hodge substructure of $T_a$,
which contradicts that $T_a$ is simple, unless it is all of $T_a$. Since the eigenvalues of $\sigma^*$ can only be $\zeta_n^k$, $k=0,\ldots,n-1$, and $n>3$ is an odd prime, this implies that $\lambda=2$ and that
$\sigma^*$  induces the identity map on $T_a$.
As $\tau^*$ is also the identity on $T_a$, we find that $D_n$ acts trivially on $T_a$. But then $T_a$ is
isomorphic to a Hodge substructure of $\tilde{\mE}/D_n$. However this is a rational surface whereas $T_a^{2,0}\neq 0$.

The Hodge endomorphism $\sigma^*+(\sigma^{-1})^*$ thus generates a subalgebra of $\End_{\Hod}(T_a)$
which is isomorphic to $F$.
In particular, $F$ acts by Hodge endomorphisms on $T_a\cong T_{\mE_a,\QQ}$. 
\qed

\subsection{Cycles inducing the real multiplication}
The real multiplication by $\zeta_m+\zeta_m^{-1}$
on $T_{\mE_a,\QQ}$ is a $\QQ$-linear endomorphism which is induced by the corresponding
endomorphism of $H^2(\tilde{\mE}_a,\QQ)$. For $k \in \{0,1,\ldots,n-1\}$
let
$$
\Gamma_k\,:=\,\{(x,\sigma^k(x))\,\in\tilde{\mE}_a\times\tilde{\mE}_a\,: x\,\in\,\tilde{\mE}_a\,\}~.
$$
be the graph of the order $n$ automorphism $\sigma\in D_n$ of $\tilde{\mE}_a$.

Let  $[\Gamma_k]\in H^{4}(\tilde{\mE}_a\times \tilde{\mE}_a,\QQ)$
be the cohomology class of the subvariety $\Gamma_k$.
Using the K\"unneth formula
$$
[\Gamma_k]\,=\,\sum_i[\Gamma_k]_{2d-i}\,\in\,
\bigoplus_i\, H^{2d-i}(\tilde{\mE}_a,\QQ)\otimes H^i(\tilde{\mE}_a,\QQ)~,
$$
and recall that with  Poincar\'e duality $H^{2d-i}(\tilde{X}',\QQ)\cong H^i(\tilde{X}',\QQ)^*$, one finds
$$
H^{2d-i}(\tilde{X}',\QQ)\otimes H^i(\tilde{X}',\QQ)\,\cong\,\End(H^i(\tilde{X}',\QQ))~.
$$
The endomorphism of $H^2(\tilde{\mE}_a,\QQ)$ defined by the K\"unneth component $[\Gamma_k]_2$ 
lies in $\End_\Hod(H^2(\tilde{\mE}_a,\QQ))$ since it has Hodge type $(2,2)$ and one has
$[\Gamma_k]_2=(\sigma^{-k})^*$,
the inverse is due to the definition of the action of the group $D_n$ on the cohomology, which is defined
by $g\cdot v:=(g^{-1})^*v$ to assure that $g\cdot (h\cdot v)=(gh)\cdot v$.

In particular, the action of $\zeta_n+\zeta_n^{-1}$ on $H^2(\tilde{\mE}_a,\QQ)$
is induced by the cycle $\Gamma_1+\Gamma_{-1}$ on $\tilde{\mE}_a\times\tilde{\mE}_a$.
This cycle induces one on $\mE_a\times\mE_a$ which defines the real multiplication on $T_{X,\QQ}$.

\subsection{Remark}

In general, if $X$ is a K3 surface with RM by $\QQ(\zeta_m+\zeta_m^{-1})$,
then there is a priori no  reason to assume that the real multiplication is induced
by a cycle with two irreducible components as we found for the $\mE_a$.
This suggests that such K3 surfaces are quite special among those with RM.
This is confirmed by the fact that
in various examples we do not find maximal families of RM K3 surfaces
(the dimension of the deformation space is given by Proposition \ref{prop:cyclo}).
Assuming the Hodge conjecture, there must be a cycle inducing the RM, but
the general member of such a maximal family  probably has a more complicated cycle
than in the case we considered here.

\subsection{Remark}

Dickson polynomials have the special feature of being permutation polynomials for certain finite fields,
only depending on the degree (but not on $a$).
The elliptic K3 surfaces considered in  Proposition \ref{prop:RM-ell}
therefore satisfy certain congruences for their point counts over finite fields;
this relates to the approach towards RM taken in \cite[Thm 1.1]{ElsenhansJ3}.


\section{Examples}\label{Dn_exa}
\label{s:ex}

\subsection{The case where $n$ is prime}
We consider elliptic surfaces with purely non-symplectic automorphisms of order $n$.
In case $n$ is a prime number we use the deformations given by Dickson polynomials as in
Proposition \ref{prop:RM-ell} to construct explicit families of elliptic K3 surfaces with RM.
A slight modification allows us to also handle the case $n=9$. In case $n=5$ we find a much larger
family using a variation of Proposition \ref{prop:RM-ell} in Section \ref{ss:5big}.

\subsection{The case $n=5$, approached via elliptic fibrations}
\label{ss:5}

Consider the rational elliptic surfaces given by the Weierstrass form
\begin{eqnarray}
\label{eq:RES-5}
S: \quad y^2 = x^3 + a_1 x + a_2, \qquad a_i\in k[t], \quad\deg(a_i)\leq i.
\end{eqnarray}
This family is 3-dimensional since the $a_i$ have $2+3=5$ coefficients but there are the scalings
$(x,y)\mapsto (\lambda^2x,\lambda^3y)$ and $t\mapsto \mu t$ to take into account.
The elliptic fibration on $S$ has a singular fibre of Kodaira type IV$^*$ at $\infty$
and generally Mordell--Weil lattice MWL $\cong A_2^\vee$
(see \cite[Table 8.2, No.\ 27]{SchuettS}).
In fact, solving for $x=$ const. such that the RHS of \eqref{eq:RES-5}
is a perfect square leads exactly to 6 sections of height $2/3$,
corresponding to the minimal vectors of $A_2^\vee$.

Base change by $t=s^5$ gives rise to a 3-dimensional family of K3 surfaces with
\begin{itemize}
\item
a non-symplectic automorphism $\sigma_0$ of order $5$,
\item
(generally) a singular fibre of type IV at $s=\infty$,
so $\NS\supset U \oplus A_2$,
\item
Mordell--Weil lattice MWL $\supseteq A_2^\vee(5)$,
the original Mordell--Weil lattice of $S$ scaled by $5$,
\end{itemize}
so $\rho\geq 6$ with very general equality by \ref{ss:CM}.
One can show that this family corresponds to the second family in
\cite[Table 2]{ArtebaniST}, listed under $S(\sigma)=H_5\oplus A_4$
(the invariant lattice under $\sigma_0^*$ acting on $H^2(X,\ZZ)$, isometric to the very general N\'eron--Severi lattice).

We can deform the above family by 
considering the K3 surface $X_a$ obtained as the base change of $S$ by $t=p_{5,a}(s)$ for $a\in\CC$.
This results in a 4-dimensional family of K3 surfaces
with the same very general N\'eron--Severi lattice (since the sections and reducible fibre deform).
By Proposition \ref{prop:RM-ell}, one has 
$\QQ(\sqrt 5)\subset \End_\Hod(T_{X_a,\QQ})$.
To see that this is an equality very generally, assume that there is a strictly larger field $F$
such that $F\subset \End_\Hod(T_{X_a,\QQ})$ very generally.
Then $m=[F:\QQ]\geq 4$,
so by \ref{ss:CM}, \ref{ss:moduli},
the 4-dimensional family would force rank$(T_X)\geq m\cdot (4+1)=20$, 
which is impossible since $\rho\geq 6$.

\subsection{The case $n=7$}
\label{ss:7}

According to \cite[\S 6]{ArtebaniST}, there are two 2-dimensional families
of K3 surfaces admitting a non-symplectic automorphism of order 7
(and this is the maximal dimension of such  families).
The dimension of the
$\QQ(\zeta_7)$-vector space $T_{X,\QQ}$ is then $l=3$
for the general $X$ in either family.
By Proposition \ref{prop:cyclo} there exist two $2l-2=4$-dimensional families
of K3 surfaces with RM by the cubic field $F=\QQ(\zeta_7+\zeta_7^{-1})$.
Using one of these families and Proposition \ref{prop:RM-ell} we find an explicit 3-dimensional family
of elliptic K3 surfaces with RM by $F$.

\subsection{Proof of Theorem \ref{thm2} (7)}
One of the families from \cite{ArtebaniST} is given by the Weierstrass forms
\begin{eqnarray}
\label{eq:7}
y^2 = x^3 + (b_1t^7+b_0)x + (c_1t^7+c_0), \qquad b_i, c_i \in\CC
\end{eqnarray}
(this is a 2-dimensional family once we account for scalings).
Generally, such a fibration has only one reducible singular fibre
(of Kodaira type III, located at $t=\infty$).
There is also a section $(x(t),y(t))$ with $x=-c_1/b_1$, of height $7/2$.
Hence a general $X$ in the family has $\Pic(X) = U \oplus K_7$, a lattice of rank $\rho=4$ in the notation of \cite{ArtebaniST},
and $T_{X} = U^2 \oplus A_6 \oplus E_8$.

We deform \eqref{eq:7} by replacing $t^7$ by $p_{7,a}(t)$, with $a\in\CC$ as in \eqref{eq:E_a},
to obtain a 3-dimensional family of K3 surfaces $X_a$ with
$F=\QQ(\zeta_7+\zeta_7^{-1})\subset \End_\Hod(T_{X_a,\QQ})$.
If $F\neq E:=\End_\Hod(T_{X_a,\QQ})$ then $E$ must be a field of degree at least $6$, and $\dim_ET_{X_a,\QQ}\leq 3$,
hence there would be at most $3-1=2$ moduli if $E$ is CM or $3-2=1$ moduli if $E$ is totally real, 
contradicting the count of 3 moduli we found.

As the singular fibre types stay the same generally and the section obviously deforms,
we infer that for the general deformation $\Pic(X_a) = U \oplus K_7$.
Thus the remaining claim  of Theorem \ref{thm2} (7) about the very general Picard number follows.
\qed

\

\subsection{The case $n=9$} \label{ss:9}
A complete classification of the K3 surfaces $X$ with a non-symplectic automorphism
$\sigma$ of order $9$ w.r.t.\ the fixed locus of $\sigma$ is  given in \cite{ArtebaniCV}.
We use deformations of a 1-dimensional family to find an explicit 2-dimensional family
with RM by the degree three totally real field  $\QQ(\zeta_9+\zeta_9^{-1})$.
 
\subsection{Proof of Theorem \ref{thm2} (9)}
A one-dimensional family (denoted by D2 in \cite{ArtebaniCV})
of elliptic K3 surfaces with very general $\rho=10$ (hence $d=\dim_\QQ(T_{X,\QQ})=12$) admitting
a purely non-symplectic automorphism of order 9 is given by
\begin{eqnarray}
\label{eq:9'}
y^2 = x^3 + bx + c_1 t^9 + c_0, \qquad b, c_1, c_0 \in\CC.
\end{eqnarray}
For very general $X$ one has $\Pic(X)= U \oplus A_2^4$,
but this becomes visible on the above fibration only indirectly, namely
through the fibre of type I$_0^*$ at $t=\infty$ and through the Mordell--Weil lattice
$\MWL(X) = D_4^\vee(3)$ which is induced from
the rational elliptic surface given by $s=t^3$ which is intermediate to the cyclic cover given by $u=t^9$.
(See \cite[Table 8.2, No.\ 9]{SchuettS} for the intermediate 
rational elliptic surface.)

We deform \eqref{eq:9'} by replacing $t^9$ by $p_{9,a}\, (a\in\CC)$.
To show that the Picard lattice is preserved by the deformation,
note that the fibre at $\infty$ is clearly preserved.
As for the Mordell--Weil lattice, it is well-known that
\[
p_{mn,a}(t) = p_{m,a^n}(p_{n,a}(t)).
\]
Presently, with $m=n=3$, this implies that also the deformation factors through a rational elliptic surface with
$\MWL=D_4^\vee$. Hence we get the same very general $\Pic$ as before.

To show that one obtains a 2-dimensional family of elliptic K3 surfaces with RM by $F=\QQ(\zeta_9+\zeta_9^{-1})$,
one modifies the proof of Proposition \ref{prop:RM-ell} by splitting $H^2(\tilde{\mE}_a,\QQ)$ into three summands
that are $D_n$-representations
on which $\sigma$ acts with eigenvalues $1$, primitive cube roots of unity and primitive nineth-roots of unity respectively. One shows that $T_a$ lies in the last summand using that the intermediate
$D_3$-cover is rational.

The family is maximal, since the Picard lattice is preserved by the deformation.
\qed

 \

\subsection{The case $n=11$}\label{ss:11}

According to \cite[\S 7]{ArtebaniST}, there are two 1-dimensional families
of K3 surfaces admitting a non-symplectic automorphism of order 11 (and this is the maximum dimension of such families). Since $[\QQ(\zeta_{11}:\QQ]=10$, the dimension of the
$\QQ(\zeta_{11})$-vector space $T_{X,\QQ}$ is $2$
for the general $X$ in either family.
By Proposition \ref{prop:cyclo} there exist two $2l-2=2$-dimensional families
of K3 surfaces with RM by $F=\QQ(\zeta_{11}+\zeta_{11}^{-1})$.
Using one of these families and Proposition \ref{prop:RM-ell} we find an explicit 2-dimensional family
of elliptic K3 surfaces with RM by $F$.

\subsection{Proof of Theorem \ref{thm2} (11)}
We consider the family of elliptic K3 surfaces given by the Weierstrass form
\begin{eqnarray}
\label{eq:11}
y^2 = x^3 + bx + (c_1t^{11}+c_0), \qquad b, c_1, c_0 \in\CC~.
\end{eqnarray}
This is a 1-dimensional family once we account for scalings and thus for the general $X$ we find
$\dim_{\QQ(\zeta_{11})}T_{X,\QQ}\geq 2$. For dimension reasons we must then have equality and
so the general $X$ has Picard number two and thus $\Pic(X) = U$ and $T_{X} = U^2 \oplus E_8^2$.
Generally, there is only one additive singular fibre (of Kodaira type II, located at $t=\infty$),
all other fibres having Kodaira type I$_1$.

We deform \eqref{eq:11}
by replacing $t^{11}$ by $p_{11,a} \, (a\in\CC)$ as in \eqref{eq:E_a}
to obtain a 2-dimensional family of K3 surfaces with RM by $F=\QQ(\zeta_{11}+\zeta_{11}^{-1})$
as before; note that this is a maximal  family, since $l=\dim_F T_{X,\QQ}=4$ and thus there are $l-2=2$ moduli.
Theorem \ref{thm2} (11) follows.
\qed

\

\subsection{Proof of Theorem \ref{thm2} (5)}\label{ss:5big}
In \cite[\S 5]{ArtebaniST} one finds a description of a family $\mA$ (case 5A)
of K3 surfaces with a non-symplectic automorphism $\sigma$ of order five.
Whereas in Proposition \ref{prop:RM-ell} we found such deformation by base change
of an elliptic fibration, we now have to consider a generalization of the proof which is based
on a $D_n$-type deformations of the quotient map $X\mapsto X/\sigma$ for a general member of this family.

The general member in the  family $\mA$ has Picard rank is two and the Picard lattice $H_5$
is generated by the classes of two smooth rational curves
(cf.\ \cite[\S 1]{ArtebaniST}).
The very general transcendental rational Hodge structure $T_X$
is thus a $\QQ(\zeta_5)$-vector space of dimension $l=(22-2)/4=5$ and the family has $l-1=4$ moduli.

The general member $\mA_p$ of the family $\mA$ is defined as the double cover of $\PP^2$,
with homogeneous coordinates $(x_0:x_1:x_2)$, branched over a smooth sextic curve $C_p$,
where $p\in\CC[x_0,x_1]$ is a degree six polynomial with six distinct zeroes in
$\PP^1_{(x_0:x_1)}$, not divisible by $x_1$:
\begin{eqnarray}
\label{eq:5-A}
C_p:\qquad p(x_0,x_1)\,+\,x_1x_2^5\,=\,0~.
\end{eqnarray}
(Changing the coordinate $x_0$, $p$ can be put in the form
$x_0(x_0-x_1)\prod_{i=1}^4(x_0-\lambda_i x_1)$,
thus making the 4 moduli apparent.
We can also observe the generators of $\Pic(X)$ as the components of the pull-back of the line $\{x_1=0\}\subset\PP^2$.)
The automorphism $\sigma$ on $\mA_p$ is induced by the automorphism
$$
\bar{\sigma}:\;\PP^2\,\longrightarrow\,\PP^2,\qquad
(x_0:x_1:x_2)\,\longmapsto\,(x_0:x_1:\zeta_5x_2)~.
$$
These surfaces were described in \cite[7.2]{GarbagnatiP} as product-quotient surfaces.
The quotient $Y=\PP^2/\langle\bar{\sigma}\rangle$
is the (singular) weighted projective space
$\PP(1,1,5)$ which is isomorphic to a cone in $\PP^6$ over
a rational normal curve of degree $5$.
After blowing up the fixed point $(0:0:1)\in\PP^2$ and the singular point of $Y\cong\PP(1,1,5)$, one obtains a
cyclic degree five cover of the Hirzebruch surfaces $\FF_1\rightarrow\FF_5$.

We deform the family $\mA$ by replacing $x_2^5$ in (\ref{eq:5-A}) by the Dickson polynomial $p_{5,a}(x_2)$
where $a\in \CC[x_0, x_1]$ is homogeneous of degree two, so that $p_{5,a}(x_2)\in\CC[x_0,x_1,x_2]$ is
homogeneous of degree $5$.
One obtains a covering map $\FF_1\rightarrow\FF_5$ with monodromy group $D_5$.
This induces a degree five covering between the double covers and similar to the proof of Proposition
\ref{prop:RM-ell} one finds that any $\mA_p$ deforms to a K3 with RM by $\QQ(\sqrt{5})$.

Since the transcendental rational Hodge structure of a general deformation still has dimension
$20$ over $\QQ$, it has dimension $l=20/2=10$ over
$\QQ(\sqrt{5})=\QQ(\zeta_5+\zeta_5^{-1})$.
Therefore the deformations of K3 surfaces in the family $\mA$ with RM by $\QQ(\sqrt{5})$
have $l-2=8$ moduli. However, the $D_n$-type deformations depend on 3 parameters, the coefficients of $a$,
so we get $4+3=7$ moduli for the $D_n$-type deformations.
\qed

\medskip
Since these K3 surfaces have Picard number $2$, they should be (very special) members of an  $8$-dimensional family
of K3 surfaces with $RM$ by $\QQ(\sqrt{5})$.
Unfortunately we do not know the general member in such an $8$-dimensional family explicitly.

\


\section{An approach using isogenies}
\label{s:isogeny}

Exploiting isogenies forms a classical approach towards exhibiting
explicit elliptic curves with CM.
More generally, it is very useful for the study of $\QQ$-curves.
We shall explore similar ideas for the elliptic fibrations over $\PP^1$ below
in order to find K3 surfaces with RM (induced by suitable rational self-maps).

\subsection{Degree 2 isogenies}
\label{ss:2-isog}

If $E$ is an elliptic curve with a 2-torsion point $P$ over a field $K$ of characteristic $\neq 2$,
then $E$ can be converted to the standard form, with $P=(0,0)$:
\begin{eqnarray}
\label{eq:E}
E: \;\;\; y^2 = x(x^2 + 2ax + b), \;\;\; a,b \in K,\;\;  b(a^2-b)\neq 0.
\end{eqnarray}
Quotienting by translation by the 2-torsion section $P$, $E$ admits a 2-isogeny (\cite[III.4.5]{Silverman1})
\[
E \stackrel{\psi}{\longrightarrow} E',\qquad
(u,v)\,=\,\left(\frac{y^2}{x^2},\,\frac{y(x^2-b)}{x^2}\right)
\]
to the elliptic curve $E'$ given by
\[
E': \;\;\; 
v^2 = u(u^2-4au+4(a^2-b)).
\]
Asking for $E$ and $E'$ to be isomorphic over $\bar K$ generally leads to the CM-curves
with j-invariants $1728, -3375$ and $8000$ (\cite[II.2]{Silverman2}.
In the realm of elliptic surfaces, i.e.\ with $K=k(t)$, 
we can set up the fibration \eqref{eq:E} to be isotrivial with general fibre one of the  above three elliptic curves $E$.
Then the surface automatically acquires CM
(and in the K3 case, at least for $j\neq 1728$,
 it turns out to be a Kummer surface for the product of $E$ with another elliptic curve).

However,
there is more flexibility since requiring that the underlying surfaces $X$ and $X'$ are isomorphic
does not imply that the isomorphism acts as identity on the base of the fibration; it need not even preserve 
the fibration.
In what follows, we  impose the condition that there is an automorphism $\sigma$ of $\PP^1$
such that the fibration $\pi$ on $X$ and the twisted fibration $\sigma\circ\pi'$ on $X'$
are isomorphic as elliptic fibrations:
\begin{eqnarray}
\label{eq:twist}
X\stackrel{\pi}{\longrightarrow} \PP^1 \;\;\; \text{ and } \;\;\;
X'\stackrel{\sigma\circ\pi'}{\longrightarrow} \PP^1~.
\end{eqnarray}
For this, note that, generally, $X$ will have singular fibres of type I$_2$ at the zeroes of $b$,
and of type I$_1$ at the  zeroes of $a^2-b$, as displayed below.

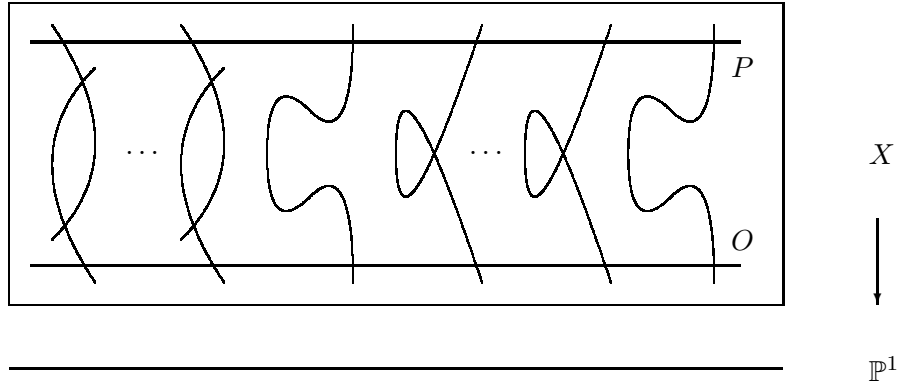
\begin{figure}[ht!]
\setlength{\unitlength}{.45in}
\begin{picture}(11,4.3)(-.5,-.1)
\thicklines

\put(0.25,3.8){\line(1,0){8.25}}
\put(8.4,3.4){$P$}

\put(0.25,1.2){\line(1,0){8.25}}

\thinlines

 \qbezier(.5,1.5)(1.5,2.5)(.5,4)
 \qbezier(1,3.5)(0,2.5)(1,1)
 
 \put(1.35,2.5){$\hdots$}
 
  \qbezier(2,1.5)(3,2.5)(2,4)
 \qbezier(2.5,3.5)(1.5,2.5)(2.5,1)


\qbezier(4.5,2.5)(4.5,1)(5.5,4)
\qbezier(4.5,2.5)(4.5,4)(5.5,1)

\qbezier(6,2.5)(6,1)(7,4)
\qbezier(6,2.5)(6,4)(7,1)

\put(5.35,2.5){$\hdots$}


\qbezier(3,2.5)(3,3.5)(3.5,3)
\qbezier(3.5,3)(4,2.5)(4,4)

\qbezier(3,2.5)(3,1.5)(3.5,2)
\qbezier(3.5,2)(4,2.5)(4,1)

\qbezier(7.2,2.5)(7.2,3.5)(7.7,3)
\qbezier(7.7,3)(8.2,2.5)(8.2,4)

\qbezier(7.2,2.5)(7.2,1.5)(7.7,2)
\qbezier(7.7,2)(8.2,2.5)(8.2,1)

\thinlines
\put(0,0.75){\framebox(9,3.5){}}

\put(0,0){\line(1,0){9}}
\put(10,0){\makebox(0,0)[l]{$\PP^1$}}
\put(10,2.5){\makebox(0,0)[l]{$X$}}
\put(10.1,1.75){\vector(0,-1){1}}

\put(8.4,1.5){\makebox(0,0)[l]{$O$}}

\end{picture}
\end{figure}

On $X'$ the singular fibres are interchanged, so we basically want to undo this using $\sigma$.
In practice, we will take $\sigma$ as an involution of $\PP^1_t$
which we 
normalize to be $\sigma(t) = -t$.
Then for \eqref{eq:twist} to give isomorphic elliptic surfaces
requires that 
\[
a = \pm a^\sigma \;\;\; \text{ and } \;\;\; b + b^\sigma = a^2.
\]

Now we restrict to the K3 setting with $k\subset\CC$ and draw the desired consequences for RM and CM.

\subsection{Proposition}
\label{prop:2}
Let $\alpha,\beta\in \CC[x]$ with $\deg(\beta)\leq 3$. Then,
\begin{itemize}
\item
if $\deg(\alpha)\leq 2$, the 5-dimensional family of elliptic K3 surfaces
$$
y^2 = x\Big(x^2 + 2\alpha(t^2)x + \frac 12 \alpha(t^2)^2 + t\beta(t^2)\Big)
$$
has CM by $\QQ(\sqrt{-2})$ and very generally the Picard rank is  $\rho= 10$;
\item
if $\deg(\alpha)\leq 1$, the 4-dimensional family of elliptic K3 surfaces
$$
y^2 = x\Big(x^2 + 2t\alpha(t^2)x + \frac 12 t^2\alpha(t^2)^2 + t\beta(t^2)\Big)
$$
has RM by $\QQ(\sqrt{2})$ and very generally the Picard rank is  $\rho= 10$.
\end{itemize}
In particular, the second point implies Theorem \ref{thm} (2).

\subsection{Remark}

The K3 surfaces in Proposition \ref{prop:2}
admit rational self-maps of degree $2$, as we will see in the proof.
They should thus be of independent interest, cf.\ \cite{Dedieu}.

\subsection{Proof of Proposition \ref{prop:2}}
\label{ss:proof-2}

The degree bounds ensure that the elliptic surfaces are K3 surfaces for general $\alpha,\beta$
(cf.\ \cite[Prop.\ 5.51]{SchuettS}).

In the first case, $a^\sigma(t)=\alpha((-t)^2)=a(t)$ and $b^\sigma(t)=\frac 12 \alpha((-t)^2)^2 + (-t)\beta((-t)^2))=
\frac 12 \alpha(t^2)^2 - \beta(t^2))$, hence $(b^\sigma+b)(t)=\alpha(t^2)^2=a^2(t)$ and we can extend $\sigma$ to an isomorphism
\begin{eqnarray*}
\varphi: \;\;\; X' \;\; & \stackrel{\cong}\longrightarrow &  \;\;\;\;\;\;X\\
(u,v,t) & \mapsto & (-2x,2\sqrt{-2}y,-t).
\end{eqnarray*}
Thus we obtain a self-map $\varphi\circ\psi$ of $X$ of degree $2$.
Since the isogeny $\psi$ preserves the regular 2-forms,
$(\varphi\circ\psi)^*$ acts on $\omega = dx \wedge dt/y$ as multiplication by   $\sqrt{-2}$.
This proves the claimed CM by $\QQ(\sqrt{-2})$.

In the second case
$a^\sigma=-a$ and $b^\sigma+b=a^2$,
so the analogous argument applies to the isomorphism
\begin{eqnarray*}
\varphi': \;\;\; X' \;\; & \stackrel{\cong}\longrightarrow &  \;\;\;\;\;\;X\\
(u,v,t) & \mapsto & (2x,2\sqrt{2}y,-t).
\end{eqnarray*}
Hence $\QQ(\sqrt 2)\subset$ End$_\text{Hod}(T_X)$ by inspection of the degree $2$ self-map $\varphi'\circ\psi$ of $X$
and its induced action on $\omega$.

We continue by verifying the stated moduli dimensions.
They amount to a simple parameter count,
compared against the 2 degrees of freedom left by scaling on the one hand $t$ and 
on the other hand admissibly $(x,y)$
(since the M\"obius transformations have to preserve the fixed points $0, \infty$ of the involution $\sigma$).
Thus the stated moduli dimensions follow.

The  bounds for the Picard numbers are an immediate application of the Shioda-Tate formula:
generally, at the zeroes of $b$, there are 8 reducible fibres of Kodaira type I$_2$ in the CM case of set-up (1),
resp.\ 6 fibres of Kodaira type I$_2$ and two fibres of type III in the RM case of set-up (2),
so in either case, we have
\[
\NS(X) \supset U \oplus A_1^8
\]
 of rank 10 at least.
But then, taking into account the moduli dimensions, End$_\text{Hod}(T_X)$ can be at most quadratic by
\eqref{eq:CM-dim}, \eqref{eq:RM-dim}.
Therefore we obtain very generally
$\rho=10$ and  CM by $\QQ(\sqrt{-2})$ resp.\ RM by $\QQ(\sqrt{2})$.
\qed

\subsection{Remark}

We emphasize that with the given Picard number (or lattice polarization),
Proposition \ref{prop:2} exhibits maximal dimensional families of K3 surfaces with RM or CM, again by \eqref{eq:CM-dim}, \eqref{eq:RM-dim}.

\subsection{Noether--Lefschetz loci}

One can easily exhibit several Noether--Lefschetz loci of the above family.
Concentrating on the RM case from Proposition \ref{prop:2},
there are three fibres of types I$_1$, I$_2$, III merging to I$_0^*$ when $t\mid \beta$ or $\deg(\beta)<3$.
This gives 3-dimensional families with very general $\rho=12$, again with RM by $\QQ(\sqrt 2)$ by construction.

Analogous results hold when we merge two pairs of I$_2$'s and I$_1$'s to I$_4$ and I$_2$
(which is easily implemented by solving for $b$ to admit a 4-fold zero at $t=1$
as the resulting equations are linear in the coefficients of $\beta$)
or when we impose additional sections (which is tedious, but doable for height $2$, for instance).

\subsection{Higher CM strata}
\label{ss:strata}

We can also find subfamilies with  CM by fields of a higher degree.
Notably, this occurs when $\alpha\equiv 0$ as then the generic fibre of \eqref{eq:E}
acquires an automorphism of order 4 and thus has CM itself.
Hence we obtain a 2-dimensional family with CM by $\QQ(\sqrt 2, \sqrt{-1}) = \QQ(\zeta_8)$
(and very general $\rho=10$ for dimension reasons).

Specializing further so that $t\beta(t^2) = p_{5,a}(t)$ resp.\ $p_{7,a}(t)$,
we obtain isolated K3 surfaces (since we can normalize $a=1$)
 with CM by $\QQ(\sqrt 2, \sqrt{-1}, \sqrt 5)$ resp.\ $\QQ(\sqrt 2, \sqrt{-1}, \zeta_7+\zeta_7^{-1})$.

Similarly, at $\alpha=1$,  $\beta=t^2$, the K3 surface admits a non-symplectic automorphism of order 3,
so it has CM by $\QQ(\sqrt 2, \sqrt{-3})$.
Since it has  singular fibres of type I$_0^*$ at $t=0$ and III$^*$ at $t=\infty$,
we conclude that $\rho=18$.

Along the same lines, for  $\alpha=1$, $\beta=t^3$, the K3 surface admits a non-symplectic automorphism of order 5,
so there is CM by $K=\QQ(\sqrt 2, \zeta_5)$.
As $\rho\geq 10$ by construction,
$T_{X,\QQ}$ can presently only have dimension $[K:\QQ]=8$ by
\ref{ss:CM} and thus $\rho=14$.


\section{Higher degree isogenies}
\label{s:higher}

It turns out that an isogeny between elliptic surfaces with torsion points of higher order that are rational over the base does not give rise to K3 surfaces with RM
(partly because those isogenies force a relatively large Picard number, whereas \eqref{eq:rho} shows that the
Picard rank is at most $16$ if a K3 has RM),
but isogenies still do the job since we only need a subgroup, the kernel of the isogeny, to be rational.
For brevity, we focus on the degree $3$ case.

\subsection{Degree 3 isogenies}

Following \cite{Top} (or \cite[II.4 \S 2]{Fricke}), one can write an elliptic curve $E$
admitting a 3-isogeny over a field $K$ of characteristic $\neq 2,3$ as
\[
E: \;\;\; y^2 = x^3 + 27a(x-4b)^2, \qquad  a,b\in K.
\]
Here the isogenous curve $E'$ is given by
\[
E':\;\;\; v^2 = u^3 -27^2 a (u-108(a+b))^2,
\]
and the 3-isogeny is
\[
(x,y) \mapsto \left(\frac 9{x^2}\left(2y^2 + 2ab^2 -  x^3 - \frac 23 ax^2\right), 27\frac{y}{x^3}(-4abx + 8ab^2 - x^3)\right).
\]
Following the approach of \ref{ss:2-isog},
we obtain the analogous cases in which $E\cong E'$
in terms of auxiliary polynomials $\alpha, \beta \in k[t]$:
\begin{enumerate}
\item
$a=\alpha(t^2), \;\; b=-\frac 12 a + t\beta(t^2)$;
\item
$a=t\alpha(t^2), \;\; b= -\frac 12 a + \beta(t^2)$.
\end{enumerate}

In the K3 setting, we derive the following families for $k\subset\CC$:

\subsection{Proposition}
\label{prop:3}
\begin{enumerate}
\item
If $\deg(\alpha)\leq 2$ and $\deg(\beta)\leq 1$, set-up (1) leads to a 3-dimensional family of K3 surfaces
with $\rho\geq 10$ and CM by $\QQ(\sqrt{-3})$;
\item
if $\deg(\alpha)\leq 1$ and $\deg(\beta)\leq 2$, set-up (2) leads to a 3-dimensional family of K3 surfaces
with $\rho\geq 10$ and RM by $\QQ(\sqrt{3})$.
\end{enumerate}
As before, the second point implies most of Theorem \ref{thm} (3).

\subsection{Proof of Proposition \ref{prop:3}}

The proof follows the same lines as \ref{ss:proof-2}.
Note that generally $X$ has
\begin{itemize}
\item
4 fibres of type II at the zeroes of $a$;
\item
4 fibres of type I$_3$ at the zeroes of $b$;
\item
4 fibres of type I$_1$ at the zeroes of $a+b$.
\end{itemize}
Hence Shioda--Tate again gives $\rho\geq 10$
since $\NS(X)\supset U\oplus A_2^4$,
and the claimed very general Hodge endomorphisms algebra
follows from the analogous parameter count using \eqref{eq:CM-dim}, \eqref{eq:RM-dim}.
\qed

\subsection{Proof of Theorem \ref{thm} (3)}
\label{ss:proof-3}

The theorem follows almost completely from Proposition \ref{prop:3}.
There is only the statement about the very general Picard number missing.
To prove this, it suffices to exhibit a special member $X$ of the family with $\rho(X)=10$.
This can be achieved by computing $\rho(X\otimes\bar\FF_p)$
at a prime $p$ of good reduction.
By the (proven) Tate conjecture, the Picard number of $X\otimes\bar\FF_p$
is encoded in the zeta function which can be computed using Magma's built in functionality for elliptic curves over function fields  \cite{Magma},
for instance.

In detail, if $\alpha = 1 + 2 t$ and $\beta=3+4t+t^2$, then the characteristic polynomial of Frobenius at $p=7$ 
on $H^2_{\text{\'et}}(X\otimes\bar\FF_p, \QQ_\ell(1)) \; (\ell\neq p)$ is
\[
(T-1)^3 (T+1) (T^2+1) (T^4+1)
\left(T^{12}+\frac{8}{7} T^{10}+\frac{6}{7} T^{8}+T^{6}+\frac{6}{7} T^{4}+\frac{8}{7} T^{2}+1\right).
\]
Since the last factor is irreducible, but not integral, it cannot be cyclotomic,
so we deduce $\rho(X\otimes\bar\FF_p) = 10$ (since we knew that $\rho\geq 10$ anyway).
This completes the proof of Theorem \ref{thm} (3).
\qed

\subsection{Remark}
The 3-dimensional family with CM by $\QQ(\sqrt{-3})$ in Proposition \ref{prop:3} (1) also fails to be maximal.
As in \ref{ss:proof-3}, this can be shown by exhibiting a special member $X$ with $\rho(X)=10$, 
so the maximal dimension of the deformation space is $(12/2)-2=4$.

Let $\alpha=3+4t+t^2$ and $\beta=1+3t$.
This leads to
the characteristic polynomial of Frobenius at $p=5$ 
on $H^2_{\text{\'et}}(X\otimes\bar\FF_p, \QQ_\ell(1)) \; (\ell\neq p)$
being
\[
(T-1)^6 (T^2+T+1)^2 \left(T^{12}-T^{10}+T^{8}-\frac{7}{5} T^{6}+T^{4}-T^{2}+1\right)~.
\]
Again we infer that $\rho(X\otimes\bar\FF_p) = 10$.

\subsection{Higher degrees}
\label{rem:5,7}
The analogous approach for isogenies of degree 5 or 7,
based on the classical Fricke parametrizations (cf.\ \cite[II.4 \S 3]{Fricke}), 
gives 1-dimensional families with CM
by $\QQ(\sqrt{-5})$ and $\QQ(\sqrt{-7})$, which we will not give here.
An isolated member of a one dimensional family with $\rho=16$ and  RM by $\QQ(\sqrt 7)$ is given by:
\begin{eqnarray*}{}
y^2&=& x^3 -\,a(t)x\,+ \,b(t),\\
a(t)&=& 27(t^2 + 13t + 49)(t^2 + 5t + 1)(t^2 - 49)^2,\\ 
b(t)&=&54(t^2 + 13t + 49)(t^4 + 14t^3 + 63t^2 + 70t - 7)(t^2 - 49)^3.
\end{eqnarray*}
This admits a self-map of degree $7$ induced by $t\mapsto 49/t$
which respects the common factors of $a$ and $b$. It acts on the regular 2-form as multiplication by $\sqrt 7$,
thus providing the RM structure

\medskip

Using 5-isogenies,
we also find the following proposition:

\subsection{Proposition}
\label{prop:CM-4}
\label{prop:RM5}

The K3 surface 
\begin{eqnarray*}
X: \quad y^2\,=\, x^{3}-27 \left(t^{2}-125\right)^{2} \left(t^{2}+10 t+5\right) \left(t^{2}+22 t+125\right) x\\
~~~-54\left(t^{2}+4 t-1\right) \left(t^{2}-125\right)^{3} \left(t^{2}+22 t+125\right)^{2}
\end{eqnarray*}
has Picard number $18$ and CM by $\QQ(\sqrt 5, \sqrt{-2})$

\subsection{Proof}
The elliptic K3 surface $X$ admits a rational self-map $g$ of degree $5$ given 
by a 5-isogeny which acts as $t\mapsto 125/t$ on the base
and as multiplication by $\sqrt 5$ on the regular 2-form.
By construction, we thus have $\QQ(\sqrt{5})\subset\End_\Hod(T_{X,\QQ})$.

The singular fibre types I$_5$, I$_1$, III twice and I$_0^*$ twice 
imply
that
\[
\NS(X) \supseteq U \oplus A_4 \oplus A_1^2 \oplus D_4^2.
\]
By the Shioda--Tate formula, this gives $\rho(X)\geq 16$
which would still be compatible with $X$ having RM.
However, the Picard number turns out to be $\rho(X)=18$, 
and consequently $X$ has CM, as evidenced by the following:

On the one hand,
the elliptic fibration admits a section of height $2$ with $x$-coordinate $-3 \left(t^{2}-125\right) \left(t^{2}+16 t+35\right)$
(and the pull-back by $g^*$ of height $10$), so $\rho(X)\geq 18$.

On the other hand,
the characteristic polynomial of Frobenius on $H^2_{\text{\'et}}(X\otimes\bar\FF_p, \QQ_\ell(1))$ at 
the ordinary prime $p=11$ 
admits the irreducible factor
\[
h = T^{4}-\frac{12}{11} T^{3}+\frac{18}{11} T^{2}-\frac{12}{11} T+1,
\]
so $\rho(X) \leq \rho(X\otimes\bar\FF_{11})\leq 18$, yielding the claimed equality.

It follows from \ref{ss:RM} that $X$ has CM by a field $F$ of degree $4$ containing $\QQ(\sqrt 5)$.
As $h$ splits completely over $\QQ(\sqrt 5, \sqrt{-2})$,
by \cite{Taelman} it can only have CM by this degree $4$ CM field.
\qed

\subsection{Remark}
As exploited in the proof of Proposition \ref{prop:CM-4},
all the  K3 surfaces from Propositions \ref{prop:3}, \ref{prop:RM5} and from \ref{rem:5,7} 
admit rational self-maps of degree 3, 5, 7 respectively
(which are not induced from the generic fibres of some isotrivial elliptic fibration).

{\renewcommand{\arraystretch}{1.3}
$$
\begin{array}{rcl}

\end{array}
$$
}

\end{document}